\documentclass[psamsfonts]{amsart}

\usepackage{amsmath,amssymb}
\usepackage[all,cmtip]{xy}

\newtheorem{thm}[subsection]{Theorem}
\newtheorem{lem}[subsection]{Lemma}
\newtheorem{prop}[subsection]{Proposition}
\newtheorem{cor}[subsection]{Corollary}

\newtheorem{rmk}[subsection]{Remark}
      {\begin{rmk}}
      {\end{rmk}}

\newcommand{\cal}{\mathcal}
\newcommand{\be}{\begin}
\newcommand{\en}{\end}

\newcommand{\mathbfr}{\mathrm}
\newcommand{\Pow}{\textsf{Pow}}

\newcommand{\CZF}{{\mathbfr{CZF}}}

\newcommand{\F}{{\cal{F}}}
\newcommand{\comment}[1]{}

\def\sub{\subseteq}

\begin{document}
\title{Constructive strong regularity and the extension property of a compactification}

\author{Giovanni Curi}
\email{giovanni.curi@unipd.it}
\date{}

\begin{abstract}
In contexts in which the principle of dependent choice may not be available, as toposes or Constructive Set Theory, standard locale theoretic results related to complete regularity may fail to hold. To resolve this difficulty, B. Banaschewski and A. Pultr introduced strongly regular locales. Unfortunately, Banaschewski and Pultr's notion relies on non-constructive set existence principles that hinder its
use in Constructive Set Theory. In this article, a fully constructive formulation of strong regularity for locales is introduced by replacing non-constructive set existence with coinductive set definitions, and exploiting the Relation Reflection Scheme. As an application, every strongly regular locale $L$ is proved to have a compact regular compactification. The construction of this compactification is then used to derive the main result of this article: a characterization of locale compactifications  (and thus, classically, of  the compactifications of a space) in terms of their ability  of extending continuous functions with compact regular codomains. Finally, an open problem related to the existence of the compact regular reflection of a locale is presented.
\end{abstract}

\keywords{Compactifications,  locales, strong regularity, inductive and coinductive definitions, dependent choice, Constructive Set Theory.}
\subjclass{06D22, 54D35, 03D70, 03F65, 03E25}

\maketitle

\section*{Introduction}
Among the different motivations for pursuing a point-free development of general topology \cite{J82}, a pivotal one is that point-free spaces (locales) often allow for choice-free proofs of theorems that in their usual topological formulation require the application of some form of choice principle. This feature makes locale theory particularly suited for contexts in which no choice principle is available, as (the internal language of) toposes \cite{J02} and constructive set theories.

This ability of locales  however fails in connection with the point-free version of a well-known standard topological result:  the proof that a compact regular locale is completely regular relies indeed on the use of the principle of Dependent Choice (DC), as for the corresponding fact for spaces \cite{J82}. As already recalled, there are contexts in which DC is not available; specifically, there are toposes, that do not validate DC, in which a compact regular locale fails to be completely regular.

For obtaining choice-free analogues of this and other results involving complete regularity, Banaschewski and Pultr \cite{BP} introduced the notion of strongly regular locale. With DC, a  locale  is strongly regular if and only if it is completely regular. Without DC, every completely regular locale is strongly regular, and  compact regular locales are strongly regular.

To formulate strong regularity, Banaschewski and Pultr use the existence of a largest interpolative relation contained in a given binary relation $R$, defined as the union of all interpolative relations contained in $R$. Unfortunately, this definition is impredicative or circular, as the defined relation is itself an interpolative relation contained in $R$. For this construction, strong set existence principles, as the Powerset axiom or
the unbounded Separation scheme, are implicitly used. While  these principles are available in toposes and classical set theory, they are regarded as highly problematic in those settings that, beside using intuitionistic logical principles (as the internal logic of toposes), also adopt a constructive notion of set.

The main set-theoretical system of this latter type is Constructive Set Theory, CST \cite{AR}\footnote{Constructive Set Theory is here used as a collective name for Myhill - Aczel systems, in particular for CZF, i.e., Constructive Zermelo Fraenkel set theory, and its standard subsystems and extensions; the specific system to be adopted is described informally in Section 1, and formally in the Appendix.}. In this article,  using coinductive methods, and exploiting the Relation Reflection Scheme in CST, I first formulate a constructive notion of strong regularity, that will allow for the choice-free proof that a compact regular locale is strongly regular, and that makes no use of non-constructive set existence principles. This appears in Section \ref{strongregularity}.
In Section \ref{compacregularcomp}, as an application of this concept, I show that, as for toposes and ZF, in CST every  strongly regular locale $L$ has a compact regular compactification, so that a locale has a compact regular compactification iff it is strongly regular   (this result is largely a reformulation of topos-valid results I presented in \cite{CuSC}, the missing ingredient there being indeed the formulation of strong regularity developed in this article).

Using the construction of this  compactification,  the core locale-theoretic (and topological) result of  this article is then derived: as is well-known, Stone-\v{C}ech compactification is characterized by its ability to allow for the (unique) extension of every continuous function with compact and regular co-domain. It is shown here that, in a similar vein,
every compactification $k$ can  be characterized  as the minimal compactification satisfying an extension property for an associated class $\cal C_k$ of continuous maps with compact regular co-domains, a fact that does not seem to have been noticed before, or adequately emphasized.
It is also proved that for two compactifications $k_1,k_2$ of a locale $L$, $k_1\leq k_2$  if and only if $\cal C_{k_1}\subseteq \cal C_{k_2}$, so that the larger the compactification, the larger the class of
continuous maps that can be extended. In particular, if $k=\beta$ is Stone-\v{C}ech compactification, $\cal C_\beta$ coincides with the class of all continuous functions on $L$ with compact and regular co-domain.\footnote{The corresponding results for compactifications of topological spaces, that follow from these in the presence of classical logic and the Prime Ideal Theorem,  do not seem to have been noticed before as well.}

In Section \ref{openproblem} an open problem is presented. In classical set theory and in toposes one can define the compact regular reflection of a locale  \cite{BM80,BP}.
A natural question is then whether this result is also derivable in CST. That this is not the case for the compact regular reflection in its full generality follows from the proof in \cite{CuExistenceSC} that
the compact \emph{completely} regular reflection of a locale $L$ (Stone-\v{C}ech
compactification of $L$) is independent from  CST, also when extended with (various principles, including) the principle of dependent choice DC. However, one can construct in CST  the compact completely regular reflection of a locale $L$ for a wide class of locales, and the locales for which this reflection is definable can  be characterized \cite{CuSC}.  So the question naturally arises whether analogous results may be obtained for the compact regular reflection. This appears to be an interesting open problem in CST.

This article has been written having in mind a reader familiar with locale theory, and possibly toposes, but not necessarily with CST.
Constructive Set Theory is a subtheory of classical set theory, so the adopted formalism and many constructions are the familiar ones of classical mathematics. CST uses  however principles (provable in classical set theory) with which the reader may not be familiar. For this reason, I included in Section \ref{preliminaries} an informal presentation of the principles of CST to be exploited, focusing on those that differ the most from the standard classical ones (the formal system for CST is recalled in Appendix). In particular, the needed principles for inductive and coinductive definitions of sets are recalled;  as examples of applications of these principles  the inductive definition of the set of finite subsets of a given set $X$, and  the coinductive definition of the largest interpolative relation $R_0$ contained in a relation $R$, are presented.

An ambition of this article is to show how, avoiding non-constructive set existence principles as Powerset and unrestricted Separation, inductive and coinductive definitions emerge as natural tools, that can be  exploited for obtaining more explicit and informative constructions.
Remarkably, principles of CST as  the axioms needed for inductive and coinductive definitions of sets, or the Relation Reflection Scheme,  that were introduced a priori, with no connection to a specific mathematical development, turn out to provide exactly what is needed for a constructive and choice-free theory of compactifications.

\section{Informal Constructive Set Theory}\label{preliminaries}
The intended audience of this article includes readers with some familiarity with locale theory and
intuitionistic logic (topos logic), but not necessarily with CST. In this section I give an informal presentation of the principles of constructive set theory to be exploited that may be unfamiliar to such audience. The last section contains a formal description of the axiom system(s) for CST. The reader is referred to \cite{AR2010} for a thorough introduction to the subject.

In the ordinary mathematical practice, one is rarely faced with the problem of showing that a certain collection of mathematical objects is a set. This is due to the fact that classical set theory, as well as the internal set theories of toposes, have strong set existence principles that trivially imply that the collections under consideration are indeed sets. For instance, given a set $X$,  one considers the discrete topology on $X$ which is given by the class $\Pow(X)$ of subsets of $X$.
There is no need of a proof in those contexts that such a class is a set, since the Powerset Axiom exactly ensures that this is the case.
The Powerset Axiom, as the unrestricted Separation axiom (see below), allow for circular non-constructive  definitions, that are rejected in contexts which adopts (intuitionistic logic and) a constructive concept of set.
In such constructive contexts, proving that a certain construction yields a set, rather than a mere class is, in general, much more demanding.

Constructive Set Theory, the setting of this article, is a subtheory of classical set theory, in particular the usual set-theoretical language and notation is adopted in this context.  Several of the CST axioms and schemes, as the Union or Pairing or Infinity axioms, are familiar ZF axioms (cf. Appendix), so  we here concentrate on those principles of CST that we shall exploit in the following, and that might be unfamiliar to the reader. Note that all such principles are provable in classical set theory.

As in classical set theory, given any formula $\phi$, one has the class $\{x: \phi(x)\}$ of elements satisfying $\phi$.
For instance, given any set or class $X$, one has the class $\Pow(X)=\{x: x\subseteq X\}$ of subsets of $X$. As already mentioned, in CST we do not have the Powerset axiom, so even when $X$ is a set, if it is not empty,  one will not be able to deduce that the class $\Pow(X)$ is a set.

\noindent Another principle of classical set theory that is not available in its full form in CST is the Separation axiom:
\medskip

for every formula $\phi$ and any set $A$, there is a set $X$ such that\smallskip

\qquad \qquad $X=\{x\in A: \phi(x)\}$.

\medskip

\noindent In CST, this principle is replaced by \emph{Restricted Separation}:
\medskip

for every restricted formula $\phi$ and any set $A$, there is a set $X$ such that\smallskip

\qquad \qquad $X=\{x\in A: \phi(x)\}$,

\smallskip

\noindent where a formula $\phi$ is \emph{restricted} if the quantifiers that occur in it (if any) are of the form $\forall x\in B$,  $\exists x\in C$, with $B,C$ sets.

\smallskip

\noindent Note that the lack of full Separation, implies that one may have in CST a subclass of a set that need not be a set, precisely subclasses of a set $A$ of the form  $C=\{x\in A: \phi(x)\}$ where $\phi$ is not restricted. Thus, the notation $U\subseteq X$ merely says that $U$ is a possibly proper subclass of $X$, while $U\in \Pow(X)$ expresses the fact that $U$ is a subset of $X$, as $\Pow(X)$ is the class of subsets of $X$.

\medskip

Another principle to be exploited that might not be familiar to the reader is the Strong Collection axiom.  Strong Collection is a strengthening of the more familiar ZF axiom of Replacement, stating that if $f:A\to C$ is a function from a set $A$ to a class $C$, the image of $A$ under $f$ is a set.

\emph{Strong Collection} states, more generally, that if $A$ is a set and $\phi(x,y)$ is a formula such that $(\forall x \in A)\exists y \phi(x, y)$, then there is a set $B$ such that
  \begin{eqnarray*}(\forall x \in  A)(\exists y \in B)\phi (x, y)\ \&\ (\forall y\in B)(\exists x \in A)\phi(x, y).
   \end{eqnarray*}

\noindent Note that from Replacement it follows that if $A$ is a set, $f:A\to C$ a function, then $f=\{(x,f(x)):x\in A\}$ is a set.
\medskip

A key role in CST, and specifically in this article,  is played by inductive and co-inductive definitions. An \emph{inductive definition}  is any class  $\Phi$ of pairs. A class $A$ is called \emph{$\Phi-$closed} if:
\smallskip

\centerline{$(a,X)\in \Phi$, and  $X\subseteq A$ implies $a\in A$.}
\medskip

\noindent
CST enjoys  the \emph{class induction property}:
\medskip

for each class of pairs $\Phi$, there is a smallest $\Phi-$closed class $I(\Phi)$; i.e., there is a
class $I(\Phi)$ such that:
\medskip

$1.$ $I(\Phi)$ is  $\Phi-$closed, and

$2.$ if the class $A$ is  $\Phi-$closed then $I(\Phi)\subseteq A$.

\medskip
\noindent
As an example, we construct inductively the class of \emph{finite subsets} of a class $C$. Following Kuratowski, we define $\Pow_{fin}(C)$ to be the class $I(\Phi_{fin})$, with $$\Phi_{fin}=\{(\emptyset,\emptyset)\}\cup \{(U\cup\{y\},\{U\}):U\in \Pow(C), y\in C\}.$$

\noindent A class $X$ is then $\Phi_{fin}-$closed if: the empty set belongs to $X$;
if a subset $U$ of $C$ belongs to $X$, the subset
$U\cup \{y\}$ also belongs to $X$, for every $y\in C$.

A class $K$ is a \emph{bound} for the inductive definition $\Phi$ if, for every $(x,X)\in \Phi$, there is a set $A\in K$ and an onto mapping $f:A\to X$. An inductive definition $\Phi$ is defined to be \emph{bounded} if:
\medskip

\noindent
1. $\{x\mid (x,X)\in \Phi\}$ is a set for every  $X$;

\noindent 2. $\Phi$ is bounded by a set.

\medskip

\noindent
The system for CST we adopt satisfies the \emph{Bounded Induction Scheme} (BIS):
\medskip

For each bounded class of pairs $\Phi$, the class $I(\Phi)$ is a set. In particular, if $\Phi$ is a set, $I(\Phi)$ is a set.
\medskip

\noindent As an example, using BIS we may prove that when $S$ is a set, $\Pow_{fin}(S)$ is a set. Indeed, $\Phi_{fin}$ is bounded by the set $K=\{\emptyset,
\{\emptyset\}\}$, and  $\{x\mid (x,X)\in \Phi_{fin}\}=\{\emptyset\}\cup\bigcup_{y\in S}\{X\cup\{y\}\}$ is a set for every subset $X$ of $S$ (since the union of set-indexed family of sets is a set, cf. Appendix).

By induction on $\Pow_{fin}(X)$, one may also prove that every $a\in \Pow_{fin}(X)$  is the range of a (not necessarily one-one) function $f:\textbf{n}\to a$, for some $n\in \mathbb{N}$, $\textbf{n}$ the set with $n$ elements: one considers the class $\{u\in \Pow_{fin}(X): (\exists n\in \mathbb{N})(\exists f)f:\textbf{n}\twoheadrightarrow u\}$, and concludes showing that this class is $\Phi_{fin}-$closed.
\medskip

Given an inductive definition $\Phi$ on a set $S$, i.e of the form $\Phi\subseteq S\times \Pow(S)$, a class $C\subseteq S$ is said to be $\Phi-$\emph{inclusive} if $C\subseteq \Gamma_\Phi (C)$, with  $$\Gamma_\Phi(C)=\{x \ \mid \ (\exists X)\ (x,X)\in \Phi \ \& \ X\subseteq C\}.$$

\noindent
CST has  \emph{the class co-induction property}:
\medskip

\noindent the class $C(\Phi)=\bigcup \{Y\in \Pow(S) \ \mid \ Y\subseteq \Gamma_\Phi (Y)\}$ is the largest $\Phi-$inclusive subclass of $S$, the class co-inductively defined by $\Phi$.
\medskip

\noindent In our system for CST the following principle is available.
\medskip

\noindent \emph{Set Coinduction Scheme} (SCS): when $S,\Phi$ are sets, $C(\Phi)$ is a set.
\medskip

As an example, we show that one can construct coinductively the largest interpolative relation contained in a given binary set-relation on a set. This construction will be exploited in the definition of strong regularity to be introduced.

\be{lem}\label{largestinterpolativesubrelation}
Let $A$ and $R$ be sets, $R$ a binary relation on $A$, $R\subseteq A\times A$. Then the largest interpolative relation $R_0$ contained in $R$ is a set.
\en{lem}

\be{proof} Let $R_0= \bigcup\{Y\in \Pow(R): (\forall x,y)[(x,y)\in Y \to (\exists z\in A)(x,z)\in Y\ \&\ (z,y)\in Y]\}$ the union of all interpolative set-relations contained in $R$. A priori, this definition only defines a class. We now prove that $R_0$ can be constructed as the set co-inductively defined by an inductive definition.  Let $\Phi$ be the inductive definition $$\{((x,z), \{(x,y),(y,z)\}): \{(x,z),(x,y),(y,z)\}\subseteq R\}.$$
Note that $\Phi$ is a set since $R$ is, so that the largest $\Phi-$inclusive class $C(\Phi)=\bigcup \{Z\in \Pow(R): Z\ \Phi\textrm{-inclusive}\}$ is a set by the SCS. We conclude showing that $R_0=C(\Phi)$: let $Y$ be an interpolative relation contained in $R$. To show that $Y\subseteq C(\Phi)$ one just need to show that $Y$ is $\Phi-$inclusive: Let $(x,y)\in Y$. Then there is $z\in A$ such that $(x,z)\in Y\ \&\ (z,y)\in Y$. Since $Y\subseteq R$, we conclude that $Y$ is $\Phi-$inclusive. To prove the converse, i.e., that $C(\Phi)\subseteq R_0$, since $C(\Phi)$ is a set, it is enough to prove that $C(\Phi)$ is interpolative. Let  $(x,y)\in C(\Phi)$. Then there is $Z\in \Pow(R)$, $Z$ $\Phi-$inclusive, with $(x,y)\in Z$. There is thus $X$ such that $((x,y),X)\in \Phi  \ \&\ X\subseteq Z$. Therefore, an element $z\in A$ exists such that $(x,z)\in R$, $(z,y)\in R$, and  $(x,z)\in Z$, $(z,y)\in Z$. Since $Z\subseteq C(\Phi)$ we conclude that $C(\Phi)$ is interpolative.
\en{proof}

\noindent Contrary to what one might expect, this result will not suffice for a constructive definition of the concept of strong regularity. Further work will be needed, in particular involving the application of a principle of CST called Relation Reflection Scheme, again provable from the ZF axioms (\cite[Theorem 2.5]{A2008}).

\medskip

\noindent Relation Reflection Scheme, RRS:
\smallskip

\begin{itemize}
\item[]
  For classes $S,R$ with
$R \subseteq S \times S$, if $a \in S$ and $\forall x \in S\ \exists y \in S\ R(x,y)$ then there is a set $S_0 \subseteq  S$ such that
      $a \in S_0$ and $\forall x \in S_0\ \exists y \in S_0\ R(x,y)$.
\end{itemize}

\noindent This scheme  has proved useful in obtaining choice-free derivations of results that had first been proved using principles of Dependent Choice (see \cite{A2008}).

\section{Locales constructively}\label{localesconstr} We assume familiarity with the theory of frames (locales) \cite{J82}. In this section we illustrate how the notions we shall need are formulated in the present constructive context.
Let $X$ be a set. As already noted, the complete lattice $\Omega(X)$ of open subsets of $X$ endowed with the discrete topology is the power class of $X$, $\Pow(X)$. This example shows that, constructively, we cannot expect a frame to be carried by a set of elements (see also \cite{CuriTarski}).  The  definition of set-generated frame/locale to be recalled  is the version of the frame/locale concept adopted in constructive set theory.

A pair of classes $(L,\leq)$, with $\leq$  a partial order on $L$, is a \emph{class locale} (or class frame) if every finite subset of $L$ has a least upper bound (meet),
every subset has a greatest lower bound (join), and if the following infinite distributive law is satisfied:

$$ x\wedge \bigvee U = \bigvee_{y\in U} (x\wedge y)$$ \noindent  for all $x\in L$, $U\in \Pow(L)$. A class-locale $L$ is said to be \emph{set-generated} by a set $B_L\subseteq L$ if, for every $x\in L$,
\begin{itemize}

\item [$i.$] the class $D_x=\{b\in B_L:b\leq x\}$ is a set,

\item [$ii.$]  $x=\bigvee D_x$.

\end{itemize}
\smallskip

\noindent $B_L$ is called a \emph{basis} of $L$. In a fully impredicative context as classical set theory or the internal language of a topos,  set-generated class-locales and ordinary locales come to the same thing.

In the following, a set-generated class-locale $(L,B_L)$ will simply be referred to as a \emph{locale} $L$, sometimes omitting the explicit mention of the basis $B_L$. When no confusion may arise, we also write $B$ for $B_L$. Note that if $(L,B_L)$ is set-generated, by Restricted Separation (Section \ref{preliminaries}) the restriction of $\leq$ to $B_L\times B_L$ is a set: $\leq \cap B_L\times B_L= \{(a,b)\in B_L\times B_L: a\in D_b\}$. Also,  if $B_L$ is a basis of $L$, and $B'_L\subseteq L$ is a set such that $B_L\subseteq B'_L$, then $B'_L$ also is a basis for $L$, as is easy to verify again using Restricted Separation.

A \emph{continuous map} of locales $f:L\to M$ is a function $f^-:B_M\to L$ (note the reverse direction) satisfying:
\smallskip

\begin{itemize}
\item [1.] $\bigvee_{a\in B_M}f^-(a)= 1=\bigvee B_L$,
\item [2.] $f^-(a) \wedge f^-(b)= \bigvee \{f^-(c): c\in B_M, c\leq  a, c\leq b\}$, for all $a,b\in B_M$,
\item [3.] $f^-(a)\leq \bigvee_{b\in U} f^-(b)$, for all $a\in B_M, U\in \Pow(B_M)$ with $a\leq \bigvee U$.
\end{itemize}
\noindent

\noindent The (in general proper) class of these maps is denoted by $Hom_{\textrm{Loc}}(L,M)$. Observe  that $Hom_{\textrm{Loc}}(L,M)$ is in a one-to-one correspondence with the collection of frame homomorphisms from $M$ to $L$, i.e., class functions from $M$ to $L$ preserving the frame structure: given $f:L\to M$, one extends $f^-:B_M\to L$ to $M$  by letting for  $a\in M$, $f^-[a]=\bigvee \{f^-(b): b\in B_M, b\leq a\} $. Trivially, for $b\in B_M$, $f^-[b]=f^-(b)$. Note in particular that $f^-$ preserves  meets and joins that exist in $B_M$. Composition $f\circ g:L\to N$ of two continuous maps $f:M\to N$ and $g:L\to M$ is defined as $(f\circ g)^-:N\to L$, $(f\circ g)^-(a)=\bigvee\{g^-(x):x\in B_M\ \& \ x\leq f^-(a)\}=g^-[f^-(a)]$.

A locale mapping $f:L\to M$ is \emph{dense} if, for $a\in M$, $f^-[a]=0$ implies $a=0$; $f$ is an \textit{embedding} (or a sublocale) if $f^-[\cdot]$ is onto, equivalently, if for every $a\in L$ there is $U\in \Pow(B_M)$ such that $\bigvee_{b\in U} f^-(b)=a$.

A locale $L$ is \emph{compact} iff every covering of $1=\bigvee B_L$ by basic elements (i.e., every $U\in \Pow(B_L)$ such that $1= \bigvee U$)  has a  finite  subcover. $L$ is \emph{regular} if, for all $a\in B_L$, $a= \bigvee \{b\in B_L: b\prec a\}$, where, for $x,y\in L$, $y\prec x\iff 1 = x\vee y^*$, with $y^*= \bigvee \{c\in B_L: c\wedge y=0\}$  the \emph{pseudocomplement} of $y$.

Note that $\{b\in B_L: b\prec a\}=\bigvee\{x\in B_L:(\forall b\in B_L)b\in D_{(x^*\vee a)}\}$ is a set, while
we could not define, as in the classical or topos-theoretic context, a locale to be regular if $x\leq\bigvee \{y\in L: y\prec x\}$, for all $x\in L$. Indeed, the class $\{y\in L: y\prec x\}$ need not be constructively a set, so that its join need not exist. It is easy however to see that the given definition is in the classical/topos theoretic context equivalent to the standard one. Similar considerations apply to the definition of compactness, and of the separation properties recalled below.

We shall need the following well-known facts.

\begin{lem}\label{denseoneone}
$i.$ If $f:L\to M$ is a dense embedding, $f^-[\cdot]:M\to L$ preserves pseudocomplements. $ii.$ If $f:L\to M$ is a dense continuous map, $L$ is compact  and  $M$ is regular, then $f^-[\cdot]:M\to L$ is one-one. $iii.$ If $h:L\to M$, $f,g:M\to N$ are continuous maps  with $h$ dense, $N$ regular, and $f\circ h=g\circ h$, then $f=g$.
\end{lem}

\be{proof} $i.$ For every $a$ in $M$, one has $f^-[a^*]=\bigvee\{f^-(x): x\in B_M\ \&\ x\wedge a=0\}\leq f^-[a]^*$;
for the converse, let $b\in M$ be such that $f^-[b]=f^-[a]^*$. By
density,  $f^-[b\wedge a]= f^-[b]\wedge f^-[a]=0$ implies $b\wedge
a=0$, so that  $b\leq a^*$, and $f^-[b]=f^-[a]^*\leq f^-[a^*]$.

$ii.$ We first prove that $f^-[\cdot]:M\to L$ is co-dense, i.e., that, for $a\in M$, $f^-[a]=1$ implies $a=1$. Indeed, $a=\bigvee\{x\in B_M:x\prec a\}$. Then, $f^-[a]=1=\bigvee\{f^-(x):x\in B_M\ \& \ x\prec a\}$. By compactness of $L$ there is  a finite subset $u$ of $B_M$ such that $f^-[\vee u]=1$ and $\vee u\prec a$, i.e., $1=\vee u^*\vee a$. Since $f^-[\vee u^*]\leq f^-[\vee u]^*=0$, by density $\vee u^*=0$, so that $a=\vee u^*\vee a=1$.

Now let $a,b\in M$ be such that  $f^-[a]=f^-[b]$. If $x\in B_M$ is such that $x\prec a$, i.e., $x^*\vee a=1_M$, then $1_L=f^-[x^*\vee a]=f^-[x^*]\vee f^-[a]=f^-[x^*]\vee f^-[b]=f^-[x^*\vee b]$. Since $f^-$ is codense, this gives $x^*\vee b=1_M$, i.e., $x\prec b$. Thus,  $x\prec a\iff x\prec b$ for all $x\in B_M$, and since by regularity for every $a\in M$, $a=\bigvee \{x\in B_M: x\prec a\}$, we conclude $a=b$.

$iii.$ We have to show that if $h^-[f^-(c)]=h^-[g^-(c)]$ for all $c\in B_N$, then $f^-=h^-$.  By regularity of $N$, for every $c\in B_N$, $c=\bigvee \{b\in B_N: b\prec c\}$.
If $b\prec c$ one has $b^*\vee c=1$, so that $g^-[b^*\vee c]=g^-[b^*]\vee g^-(c)=1$. Then,
$$f^-(b)=f^-(b)\wedge (g^-[b^*]\vee g^-(c))=(f^-(b)\wedge g^-[b^*])\vee (f^-(b)\wedge g^-(c)).$$
\noindent On the other hand, since $f\circ h=g\circ h$,
\smallskip

\noindent $h^-[f^-(b)\wedge g^-[b^*]]=h^-[f^-(b)]\wedge h^-[g^-[b^*]]= h^-[f^-(b)]\wedge h^-[f^-[b^*]]=h^-[f^-[b\wedge b^*]]=0,$
\smallskip

\noindent and since $h$ is dense, we get $f^-(b)\wedge g^-[b^*]=0$. Together with the above, this gives
$f^-(b)=f^-(b)\wedge g^-(c)$, whence $f^-(b)\leq [g^-(c)]$ for all $b\prec c$, so that  $f^-(c)\leq g^-(c)$. Symmetrically $g^-(c)\leq f^-(c)$.
\end{proof}

Compared with regularity, a constructive concept of complete regularity requires some more work.
The standard definition of complete regularity for locales demands  that for every $x$ in $L$,  $x\leq\bigvee \{y\in L: y\prec\prec x\}$, where $y\prec\prec x$, $y$ really inside (or completely below) $x$, if there is a scale from $y$ to $x$. A map $s:\mathbb{I}\to L$ from the rational unit interval $\mathbb{I}$  to $L$ is called a \emph{scale}  from $y$ to $x$ if $s(0)=y, s(1)=x$ and, for $p< q$, $s(p)\prec s(q)$.

\noindent In this case we have the additional difficulty that, even restricting to basic elements as we did in redefining regularity, we still not have that $\{b\in B: b\prec\prec a\}$ is a set, because, by contrast with $\prec$, the really inside relation $\prec\prec$ is not a set when restricted to $B\times B$ (as the class of scales is not a set). This difficulty is circumvented as follows.

\noindent A locale $L$ will be defined \emph{completely regular} if a family   $ri:B\to \Pow(B)$ of subsets of $B$ exists such that \emph{i.} for all $a\in B$, $a= \bigvee ri(a)$, and \emph{ii.} for all $b\in ri(a)$ a scale exists from $b$ to $a$. The idea is that, given a basic element $a$, one does not try to construct the set of all basic elements $b$ for which a scale exists from $b$ to $a$; to be able to establish that $L$ completely regular, it will suffice to have a set of such elements that is sufficiently rich to cover $a$.

As said for regularity, it is not difficult to verify that this concept is classically (or topos-theoretically) equivalent to the standard one, as well as that this notion allows for the derivation of standard desirable properties of complete regularity: it is hereditary, invariant under isomorphisms, and allows to derive Tychonoff embedding theorem \cite{CuImLoc}.

The given constructive definition of complete regularity circumvents the use of the really inside relation. However, this relation enjoys the important interpolation property that is in particular useful in the construction of completely regular compactifications.
It is possible to reconstruct a satisfactory set-analogue of the really inside relation, satisfying  the interpolation property, adding enough scales to the basis as follows.

For $a,b\in B$, let $b\prec\prec_B a$ iff a scale of basic elements exists from $b$ to $a$. Note that this relation is interpolative on $B$. We say that $L$  has a {\it completely regular basis} $B$ if, for all $a\in B$, $a= \bigvee \{b\in B: b\prec\prec_B a\}$.

Of course, a locale with a completely regular basis is completely regular. For the converse the following holds true.

\begin{lem}\label{}
If $(L, B)$ is a completely regular locale, then $L$ has a completely regular basis.
\end{lem}

\be{proof} By  definition of complete regularity, for every pair $a,b$ such that $b\in ri(a)$ a scale of elements of $L$ exists from $b$ to $a$. We may then apply Strong Collection (Section \ref{preliminaries}) to collect a set $V'$ of scales so that for all $a,b\in B$ such that $b\in ri(a)$ there is a scale $s:\mathbb{I}\to L$ in $V'$ from $b$ to $a$, and conversely for all  $s:\mathbb{I}\to L$ in $V'$ there are $a,b\in B$ such that $s:\mathbb{I}\to L$ is a scale from $b$ to $a$. Let $V$ be the set of elements of $L$ that belong to a scale in $V'$, $V=\bigcup_{s\in V'} \textrm{Range}(s)$. Since $V$ is a set, we can construct a new basis $B'=B\cup V$ for $L$. To prove that $B'$ is a completely regular base, i.e., that for all $u\in B'$, $u= \bigvee \{v\in B': u\prec\prec_{B'} u\}$,
note that we have, for all $u\in B'$, $u=\bigvee\{x\in B: (\exists b\in B) b\in D_u \ \&\ x\in ri(b)\}$.  Since for every $c\in ri(b)$ there is in $V'$ a scale from $c$ to $b$, we have $u=\bigvee\{c\in B: (\exists b\in B) b\leq u \ \&\ c\prec\prec_{B'} b\}$, and a fortiori $u=\bigvee\{v\in B': (\exists b\in B)b\leq u\ \&\ v\prec\prec_{B'} b\}$. But if $v\prec\prec_{B'} b$ for $b\leq u$, then also $v\prec\prec_{B'} u$: given the scale $s:\mathbb{I}\to L$ from $v$ to $b$, replace $s(1)=b$ with $s(1)=u$ to get a scale from $v$ to $u$.
\end{proof}
\medskip

\noindent An important property of completely regular bases is that any extension of a completely regular basis is again a completely regular basis.

\section{Strong Regularity}\label{strongregularity}

As is well-known \cite{J82},  in a compact regular locale the well-inside relation interpolates. Then, if the principle of countable dependent choice (DC) is available, whenever $x\prec y$,  a scale from $x$ to $y$ can be constructed\footnote{This holds true for systems as ZF or the internal logic of toposes. In CST, due to the absence of Powerset, a stronger version of the axiom of Dependent Choice, the axiom of Relativized Dependent Choice, RDC (see e.g. \cite{AR2010} for this principle) is needed. If the locale is compact and regular, DC suffices for the construction of a scale, since  Lemma \ref{lemmacr} below ensures that the element $z$ such that $x\prec z\prec y$ belongs to a set (a base), rather than a class.\label{footnoteRDC}}. Thus, assuming DC,  a compact regular locale is completely regular. In \cite{BP} Banaschewski and Pultr note that this result may fail to hold in toposes that do not validate DC.  For resolving this and related difficulties, they introduce strongly regular locales.

Banaschewski and Pultr begin by noting that  any binary relation $R$ contains a largest interpolative relation. Indeed, as the union of any set of interpolative relations is an interpolative relation, defining $R_0$ to be the union of all interpolative relations contained in $R$ will yield the largest interpolative relation contained in $R$. Taking $R=\prec$ one thus obtains the largest interpolative relation $\prec_0$ contained in the well-inside relation. A locale $L$ is then defined strongly regular if $x\leq\bigvee \{y\in L: y\prec_0 x\}$, for all $x\in L$.

Unfortunately, this definition, albeit adequate to the choice-free context of toposes, is impredicative or circular: $R_0$ is indeed itself an interpolative relation contained in $R$, so that it belongs to the collection of interpolative relations of which it is the union.
In the present constructive setting, defining $R_0=\bigcup\{R'\in \Pow(R): R' \ \textrm{interpolative}\}=\{(x,y)\in R: \exists R'\in \Pow(R) [R' \ \textrm{interpolative}\ \&\ (x,y)\in R']\}$ a priori only yields a proper class, since $\Pow(R)$ cannot be proved to be a set in CST when $R$ is not empty \cite{CuExistenceSC}, and the unbounded separation scheme is not available. Furthermore, in this context, as already noted, $R=\prec\subseteq L\times L$ is a proper class, rather than a set.

Exploiting the Set Coinduction Scheme  and the Relation Reflection Scheme (Section \ref{preliminaries}), in this section I present an alternative method for defining strong regularity, that is both non-circular and adequate to a choice-free intuitionistic setting.

On a locale $L$, let the class
$$\prec_0=\bigcup\{R\in \Pow(\prec): (\forall x,y) [R(x,y) \to (\exists z\in L) R(x,z)\ \&\ R(z,y)]\}$$ be the union of all interpolative set-relations contained in the well-inside relation $\prec$ on $L$.
Note again that $\prec_0$ is a class, so that it need not belong to $\Pow(\prec)$.
Similarly to what has been said before for regularity and complete regularity, one cannot define $L$ to be strongly regular if $x\leq\bigvee \{y\in L: y\prec_0 x\}$ for all $x\in L$, since the class $\{y\in L: y\prec_0 x\}$ is not a set in general. Again, we circumvent the problem by defining  $L$ to be \emph{strongly regular} if a family $si:B\to \Pow(B)$ of subsets of $B$ exists such that
\medskip

\emph{i.} for all $a\in B$, $a= \bigvee si(a)$, and

\emph{ii.} if $b\in si(a)$, then  $b\prec_0 a$.
\medskip

\noindent Trivially, in the classical and topos-theoretic context, $L$ is strongly regular according to the Banaschewski and Pultr notion if and only if it is strongly regular for this definition. Also, one has that a completely regular locale in the sense of the previous section is strongly regular (just let $si(a)=ri(a)$ for all $a\in B$). Conversely, assuming the axiom of (relativized) dependent choice, one proves that a strongly regular locale is completely regular, with $ri(a)=si(a)$ for all $a\in B$ (cf. note \ref{footnoteRDC}).

As for complete regularity, this definition is however not completely satisfactory. It avoids the circularity in the definition of the largest interpolative set-relation contained in $\prec$, but the important interpolation property of this relation is lost. To recover this property we shall argue as follows. First we construct a larger basis $B'$ containing `enough interpolants'. Over this new basis, we then construct coinductively the largest set-relation contained in the well-inside relation. This construction is the analogue for strong regularity of the construction of completely regular bases for completely regular locales.

Given a strongly regular locale $(L,B)$, with  $si:B\to \Pow(B)$ as above, let $S$ be the class  $$\{\{x_0,...,x_n\}: n\in\mathbb{N}_{>0} \ \&\ x_0,...,x_n \in L\ \&\ (\forall i<n) x_i\prec_0 x_{i+1}\}.$$
Let $R\subseteq S\times S$ be defined by $$R(\{x_0,...,x_n\}, \{y_0,...,y_m\})$$
$$\iff$$

\noindent $x_0=y_0\ \&\ x_n=y_m\ \&\ \{x_0,...,x_n\}\subseteq \{y_0,...,y_m\}\ \&\ $

\begin{flushright}
$\ \&\ (\forall k<n)(\exists k'<m)\ x_k\prec_0 y_{k'}\prec_0 x_{k+1}.$
\end{flushright}

\noindent We then have $(\forall l_1\in S)(\exists l_2\in S) R(l_1,l_2)$, since the class $\prec_0$ is interpolative. Given a pair $(c,b)\in S$, i.e., such that $c\prec_0 b$, using the Relation Reflection Scheme (Section \ref{preliminaries}) one has a set $S_0$ such that $S_0\subseteq S$, and $$(c,b)\in S_0\ \&\  (\forall l_1\in S_0)(\exists l_2\in S_0) R(l_1,l_2).$$
Thus, for all pairs $(c,b)$ such that $c\in si(b)$ there is such a set $S_0$. Let $SI=\{(c,b):c\in si(b)\}$.  By Strong Collection (Section \ref{preliminaries}), a set $V'$ exists such that
 $$(\dagger)\ \ \  (\forall (c,b)\in SI)(\exists S_0\in V')(c,b)\in S_0\subseteq S \ \&\  (\forall l_1\in S)(\exists l_2\in S) R(l_1,l_2)$$
and
$$ (\forall S_0\in V')(\exists (c,b)\in SI)(c,b)\in S_0\subseteq S \ \&\  (\forall l_1\in S)(\exists l_2\in S) R(l_1,l_2).$$
\smallskip

\noindent Let $V= \{z\in L:(\exists S_0\in V')(\exists l\in S_0)z\in l\}$.
Define $B'=B\cup V$. Note that $B'$ is a set, and a basis of $L$.

\smallskip

Let then $\prec'$ be the restriction of the well-inside relation to $B'\times B'$, $\prec' = \{(y,x)\in B'\times B': 1\leq x\vee y^*\}$.
The relation $\prec'$ is a set:  one has $1\leq x\vee y^*\iff (\forall b\in B)b\in D_{(x\vee y^*)}$. Since  $D_{(x\vee y^*)}=\{c\in B: c\leq x\vee y^*\}$ is a set, as $L$ is set-generated by $B$, by Restricted Separation (Section \ref{preliminaries}) we can conclude that $\prec'$ is a set. By  Lemma \ref{largestinterpolativesubrelation}, we can then construct the largest interpolative set-relation $\prec_0'$ contained in $\prec'$.

\begin{lem}\label{}
Let $u,v,x,y\in B'$. If $x\leq u\prec_0' v\leq y$, then  $x\prec_0' y$.
\end{lem}
\be{proof} Consider the relation $R=\{(x,y)\in B'\times B': (\exists (u,v)\in B'\times B')x\leq u\prec_0' v\leq y\}=\{(x,y)\in B'\times B': (\exists (u,v)\in B'\times B') x\in D_u\ \&\ (u,v)\in \prec_0' \ \&\ v\in D_y\}$. By Restricted Separation $R$ is a set. Moreover, $R$ is interpolative and  $R\subseteq \prec'$. Since $\prec_0'$ is the largest relation with these properties, $R=\prec_0'$.
\en{proof}

We shall say that a locale $L$ has a \emph{strongly regular basis} $B$ if the largest interpolative relation $\prec_0'$ contained in the restriction of the well-inside relation to $B$ is compatible, i.e., if, for all $a\in B$, $a= \bigvee \{b\in B: b\prec_0' a\}$.  A locale with a strongly regular basis is strongly regular. Conversely, the following obtains.

\begin{prop}\label{stronglyregularbase}
If $(L, B)$ is a strongly regular locale, then $L$ has a strongly regular basis $B'$ extending $B$.
\end{prop}

\be{proof} We prove that the basis $B'=B\cup V$ we have constructed is a strongly regular basis for $L$. We have to show that for all $u\in B'=B\cup V$, $u= \bigvee \{v\in B': v\prec_0' u\}$. Let $u\in B'$. Then $u=b\in B$, or $u\in V$.

In the first case, as $L$ is strongly regular, $b\leq \bigvee si(b)$. We show that $si(b)\subseteq \{v\in B': v\prec_0' u\}$. If $c\in si(b)$ by $\dagger$ there is $S_0\in V'$ such that
$$(*)\ \ \ (c,b)\in S_0\ \&\  (\forall l_1\in S_0)(\exists l_2\in S_0) R(l_1,l_2).$$
On $B'$ consider the binary relation
$$R'(u,v)\iff (\exists \{x_0,...,x_n\}\in S_0)(\exists i)(\exists j>i) x_i=u \ \&\  x_j=v.$$
Note that $R'(u,v)$ implies $u\prec' v$, and that $R'$ interpolates: assume $R'(u,v)$, then, by $(*)$, there is $y$ such that $R'(u,y)$ and $R'(y,v)$.
As a consequence, $R'\subseteq \prec_0'$, and since in particular  $R'(c,b)$, we conclude $si(b)\subseteq \{v\in B': v\prec_0' b\}$, whence $b= \bigvee \{v\in B': v\prec_0' b\}$.

If, on the other hand, $u\in V$, then  $u=\bigvee\{b\in B: b\leq u\}$. For every $b\in B$ such that $b\leq u$, one has $b\leq \bigvee \{v\in B': v\prec_0' b\}$. Since then $u\leq \bigvee\{v\in B':(\exists b\in B)b\in D_u\ \&\ v\prec_0' b\}$, we conclude observing that by the previous lemma $v\prec_0' b\leq u$ implies $v\prec_0'  u$.
\en{proof}

\noindent We leave to the reader to verify that every basis extending a strongly regular base is strongly regular.

This section ends with a proof that a compact regular locale is strongly regular. This will be the consequence of the following lemma, to be often invoked in Section \ref{compacregularcomp}.
Given a locale $(L, B)$, we can always extend its basis $B$ to a basis $B^*$ that is (a set and) a sub-pseudocomplemented distributive lattice of $L$, as follows.
$B^*$ is defined inductively by the following clauses:
\smallskip

\begin{itemize}
\item[]  if $b\in B$ then $b\in B^*$;

\item[] if $u\in B^*$ then $u^*\in B^*$;

\item[] if $V$ is a finite subset of $B^*$, then $\wedge V,\ \vee V\in
    B^*$,
\end{itemize}
\smallskip

\noindent where  $\wedge,\vee,$ and $^*$ are the operations in $L$. In other terms, $B^*$ is generated by the inductive definition

\begin{align*}
  \Phi_{B^*} &= \{(b, \emptyset): b\in B\}\cup\{(u^*, \{u\}): u\in L\}\cup\\
   &\cup \{(\wedge V, V): V\in \Pow_{fin}(L)\}\cup \{(\vee V, V): V\in \Pow_{fin}(L)\}
\end{align*}

\noindent Since this inductive definition is bounded, $B^*=I(\Phi_{B^*})$ is a set (Section \ref{preliminaries}).

The following lemma is proved in \cite{CuSCAlexandroff}. The proof is reformulated here for the reader's convenience.

\begin{lem} \label{lemmacr}Let $L$ be a compact regular locale, $P$ a set, a sub-lattice and a basis of $L$. Assume, for $b\in L$, $U\in\Pow(P)$, that
$b\prec \bigvee U$. Then $b=0$, or finite subsets $\{p_1,...,p_n\}$, $\{p'_1,...,p'_n\}$, and $\{p''_1,...,p''_n\}$ of $P$ exist with $p_i\prec
p'_i\prec p''_i$, $p''_i\in U$, and such that $$b\leq \vee p_i \prec \vee p'_i\prec \bigvee U.$$
\end{lem}

\be{proof} Since $L$ is regular, and since $P$ is a basis,  by $b\prec \bigvee U$, one obtains
$1=b^*\vee \bigvee \{p\in P : (\exists p''\in U)(\exists p'\in P)p'\prec p'',\ p\prec p'\}$.
By compactness, thus, $1= \vee s_0$, with $s_0\in \Pow_{fin}(P)$  a finite subset of $\{b^*\}\cup \{p\in P : (\exists p''\in U)(\exists p'\in P)p'\prec p'',\ p\prec p'\}$. By induction on $\Pow_{fin}(P)$ (Section \ref{preliminaries}),    $s_0= u_0\cup v_0$, with $u_0\subseteq \{b^*\}, v_0\subseteq \{p\in P : (\exists p''\in U)(\exists p'\in P)p'\prec p'',\ p\prec p'\}$,
and both finite.
Then $b\leq \vee v_0$. Again by induction, $v_0$ is either empty or inhabited. In the first case $b=0$; in the second, let $v_0=\{p_1,...,p_n\}$ be an enumeration of $v_0$ (Section \ref{preliminaries}).
Thus, $b\leq \vee \{p_1,...,p_n\}$. For all $i\leq n$, there are $y\in P$ and $z\in U$ such that $p_i\prec y\prec z$. Using induction on the length $n$ of the enumeration of $v_0$ (see also \cite[Lemma 8.2.8]{AR2010}, provable choice for finite sets), one constructs the required subsets $\{p'_1,...,p'_n\}$,  $\{p''_1,...,p''_n\}$.
\en{proof}

In particular, if for $a,b\in P$ one has $b\prec a$, then there are
$p,q\in P$ with $b\leq p \prec q\prec a$. Therefore, a compact regular locale $L$ is strongly regular, with $P$ a strongly regular basis for $L$.

\section{Compact regular compactifications}\label{compacregularcomp}
In contexts as toposes or constructive set theory, in which the principle of Dependent Choice is not assumed, compact regular compactifications do not coincide with compact completely regular compactifications. In this section we deal with compactifications of the first type. Part of the material in this section, formulated for the topos-theoretic context, is presented in \cite{CuSCAlexandroff}. It is here reproduced in a version adequate to CST.

A compact regular compactification of a locale $L$, in the following simply a \textit{compactification}, is a dense embedding $k:L\to M$, with compact regular co-domain.
Exploiting the coinductive formulation of strong regularity derived in the previous section, we shall show how, given a compactifiable  locale $L$, and any set-indexed family of locale mappings $\F\equiv
\{f_i\}_{i\in I}$, $f_i:L\to L_i$, with compact and regular co-domain,  a compactification $\gamma L$ of $L$ can be constructed with the property that each function in $\F$ has a unique extension to $\gamma L$.

The construction of the compactification to be presented will then allow us to show that any compactification can be
characterized in terms of its ability of extending  continuous functions with compact regular codomain,  an aspect of compactifications that does not seem to have been emphasized before. The relation of these results with the compact regular reflection of $L$ (the `weak' Stone-\v{C}ech compactification) is discussed in the next section.
\smallskip

As a compact regular locale is strongly regular, and since, as is easy to verify, strong regularity is hereditary (if $(L,B_s)$ is a strongly regular locale with strongly regular basis $B_s$, and $f:L'\to L$ is an embedding, $B'_s=\{f^-(b):b\in B_s\}$ is a strongly regular basis for $L'$), a compactifiable locale must be strongly regular. We now show that the converse also holds constructively.

In \cite{B90} Banaschewski showed that relations of strong inclusion on a frame $L$ may
be used for defining the compactifications of that frame. The whole frame structure is not needed for carrying out Banaschewski's construction. In particular, strong inclusions can be considered more generally over pcd-lattices \cite{CuSCAlexandroff}. This fact is particularly relevant in a constructive context as the present one, since, while in such a context a (non-trivial) frame is carried by a proper class,  a pcd-lattice may well  be  carried by a set. Recall in particular from the previous section that every set $B$ generating a frame can be extended to a pcd-lattice $P=B^*$ that is a set; recall also that, as noted earlier,  the partial order restricted to any set of generators is a set.

In the following, by \emph{pcd-lattice} we mean a partial order $(P,\leq)$ with $P,\leq$ sets, satisfying the standard axioms of a pseudo-complemented distributive lattice. A \textit{strong inclusion} on a pcd-lattice $P$ is a binary set-relation $\lhd$ on $P$ satisfying:
\smallskip

1. $0\lhd 0$, $1\lhd 1$;

2. if $x\leq a \lhd b\leq y$ then $x\lhd y$;

3. if $x\lhd a$, $x\lhd b$ then $x\lhd a\wedge b$;

4. if $x\lhd a$, $y\lhd a$ then $x\vee y \lhd a$;

5. if $a \lhd b$ then $b^*\lhd a^*$;

6. $\lhd \sub \prec$;

7. if $x\lhd y$ then $x\lhd z\lhd y$, for some $z\in P$.
\smallskip

\noindent Note that from 6. it follows $\lhd \subseteq \leq$. On every pcd-lattice $P$ the well inside relation $\prec$ (is a set and) satisfies conditions $1$ to $6$, while  $\prec_0'$ is a strong inclusion: that $\prec_0'$ is a set is due to Proposition \ref{largestinterpolativesubrelation}; that it is a strong inclusion is proven as
in  \cite{BP}:
1. consider $R = \{(0,0), (1,1)\}$. $R$ is interpolative and contained in $\prec$. Since $\prec_0'$ is the largest such relation, $\{(0,0), (1,1)\}\subseteq \prec_0'$.
2. Let $R= \{(x,y)\in P\times P: (\exists a\in P)(\exists b\in P)x\leq a \prec_0' b\leq y\}$. $R$ is a set, interpolative and contained in $\prec$, so $R\subseteq \prec_0'$. Properties from 3 to 5 are proven similarly. 6. and 7. hold by definition.

If $P$ is a sub-pcd-lattice of a frame $L$, and $L$ is set-generated by $P$, a strong inclusion $\lhd$ on $P$ is said to be \textit{compatible with $L$} if $a =\bigvee_L \{x\in P: x\lhd a\}$, for all $a\in P$. The following lemma provides a method for constructing strong inclusions on $P$.

\begin{lem}\label{idsi} Let $P$ be a pcd-lattice, and $R$ be an interpolative binary set-relation on $P$ with $R\sub \prec$. The
relation  $\lhd_R$ defined inductively on $P$ as the least relation containing $R$ and closed under conditions 1-5 is a strong inclusion on $P$.
\end{lem}

\be{proof} Let $\lhd_R=I(\Phi_R)$, with
\begin{align*}
  \Phi_R &= \{((x,y), \emptyset): (x,y)\in R\}\cup\\
   &\cup \{((0,0), \emptyset)\}\cup \{((1,1), \emptyset)\}\cup\\
   &\cup \{((x,y), \{(a,b)\}): x,y,a,b\in P,\ x\leq a,\ b\leq y\}\cup\\
   &\cup \{((x,a\wedge b), \{(x,a),(x,b)\}): x,a,b\in P\}\cup\\
   &\cup \{((x\vee y,a), \{(x,a),(y,a)\}): x,y,a\in P\}\cup\\
   &\cup \{((b^*,a^*), \{(a,b)\}): a,b\in P\}.
\end{align*}

\noindent
This inductive definition is a set, and so is bounded. Thus, $\lhd_R=I(\Phi_R)$ is a set.
Obviously,   $\lhd_R$ satisfies conditions from 1  to 5. Indeed, assume for instance $x\leq a \lhd_R b\leq y$ as in condition 2.
Then, since $\lhd_R$ is $\Phi_R$-closed, also $x\lhd_R y$.

The proof of $6$ and $7$  is by induction. So, to prove 6, one shows that $\prec$  is $\Phi_R$-closed. Since $\lhd_R$ is the least $\Phi_R$-closed class, we conclude $\lhd_R\subseteq \prec$: the base cases ($(x,y)\in R$ and  1) are trivial (as $R\subseteq \prec$ by hypothesis); easy calculations show that $\prec$ is closed for the other conditions in $\Phi_R$.

Proof of 7: we prove that the set $K=\{(x,y)\in \lhd_R: (\exists z\in P)(x,z)\in \lhd_R\ \&\ (z,y)\in \lhd_R \}$ is $\Phi_R$-closed, so that it must contain $\lhd_R$. The base cases are immediate. Now assume $x\leq a, b\leq y$ for $x,y,a,b\in P$, and let $(a,b)\in K$. Then, by the proof that condition 2 on strong inclusions  holds, $(x,y)\in \lhd_R$ and  $x\lhd_R z\lhd_R y$, for some $z\in P$, so $(x,y)\in K$, as wished. The closure for the other conditions in $\Phi_R$ is proved similarly, as a further example we show that if $(x,a),(x,b)\in K$ then also $(x,a\wedge b)\in K$: from $(x,a),(x,b)\in K$ we have $x\lhd_R c_1\lhd_R a, x\lhd_R c_2\lhd_R b$ for some $c_1,c_2\in P$. Then, by 3 and
2, $x\lhd_R c_1\wedge c_2\lhd_R a\wedge b$.
\en{proof}

\noindent Note that, for $R=\emptyset$, one gets the least strong inclusion on $P$: $x\lhd_\emptyset y$ iff $x=0$ or $y=1$. Recall that an ideal in a join-semilattice is a set $I$ such that
\medskip

$i.$ if $a\in I, b\leq a$, then $b\in I$,

$ii.$ if $u\subseteq I$ is finite, then $\vee u\in I$.
\medskip

\noindent A \emph{round ideal} of a pcd-lattice $P$ with a strong inclusion $\lhd$ is an ideal $I$ such that
for all $b\in I$ there is $a\in I$ with $b\lhd a$.

\noindent Ordered by inclusion,  the class  ${\cal R}(P)$ (or ${\cal R}(P,\lhd)$ when
confusion may arise) of such ideals is a class-frame, in fact a subframe of the frame of ideals of $P$. Indeed, the meet of a finite set $F$ of round ideals  $\wedge F=\cap  F$ is round, as is the join $\bigvee \cal U=\{x\in P:(\exists u\in \Pow_{fin}(P))(\forall y\in u)(\exists I\in \cal U)y\in I\ \&\ x\leq \vee u\}$, for  $\cal U\subseteq {\cal R}(P)$ a set of round ideals.

\noindent Observe also that for all $a\in P$, $\Downarrow a\equiv \{b\in P: b\lhd a\}$ is a round ideal, and that the set $B_d=\{\Downarrow a : a\in P\}$ is a basis of ${\cal R}(P)$. With this basis, $({\cal R}(P),B_d)$ is a set-generated frame.

In \cite{B90} Banaschewski proved that given a compatible strong inclusion on a frame $L$, the locale map $h:L\to {\cal R}(L)$, with $h^-(I)=
\bigvee I$, is a compactification of $L$. The round ideal construction cannot be performed over a non-trivial frame constructively (since a frame is carried by a proper class, see also Section \ref{openproblem}). Fortunately, the whole frame structure  is not needed for obtaining a compact  regular locale.

\begin{lem}\label{representation}
Given a pcd-lattice $P$ and a strong inclusion $\lhd$ on $P$, $({\cal R}(P),B_d)$ is a compact regular locale.
\end{lem}

\be{proof} Compactness is trivial.
Regularity: with minor adjustments, the proof in \cite[Lemma 2]{B90} is valid in CST if one replaces $L$ with $P$. We reproduce it here for the reader's convenience.
Note first that, if $x\lhd a$, then $\Downarrow x\prec \Downarrow a$ in ${\cal R}(P)$: indeed if $x\lhd a$ there are $b,c$  such that $x\lhd c \lhd b \lhd a$. By 5.,  also $a^*\lhd b^* \lhd c^* \lhd x^*$. Moreover, as $b\vee c^*=1$ ($1$ the greatest element in $P$), $P=\Downarrow a\vee \Downarrow (x^*)\leq \Downarrow a\vee (\Downarrow x)^*=P$ (the latter since trivially $\Downarrow (x^*)\cap \Downarrow (x)=0$).
By 7. and 4. we have  $\Downarrow a=\bigvee \{\Downarrow x:x\lhd a\}$, so that, in conclusion,  $\Downarrow a=\bigvee \{\Downarrow x\in B_d:\Downarrow x\prec  \Downarrow a\}$, proving the regularity of ${\cal R}(P)$.
\en{proof}

\medskip

\noindent Given  a compact  regular frame $L$, set-generated by a basis $B$, the frame ${\cal R}(B^*)$, obtained considering on $B^*$  the well inside relation as strong inclusion (recall Lemma  \ref{lemmacr}), is isomorphic to $L$ via  $j^-:B_d\to L$, $j^-(\Downarrow b)=\bigvee \Downarrow b$: we leave to the reader to prove that $j$ is a locale mapping; then, $j^-(\Downarrow b)=b$ by strong regularity, so that $j^-$ is onto. We want to show that $j^-[]:{\cal R}(B^*)\to L$ is \emph{one-one}; note  that $j^-[I]=\bigvee I$, so we have to prove that for all $I,J\in {\cal R}(B^*)$, if $\bigvee I=\bigvee J$ then $I=J$. This follows from the observation that if $L$ is compact, then a strongly regular ideal $I$ on $B^*$ is such that $x\in I\iff 1_L= x^*\vee \bigvee I$ (cf \cite[Lemma 5]{BM80}).

Thus, by the above lemma, compact regular frames are precisely the frames of round ideals over some pcd-lattice endowed with a strong
inclusion.

\medskip

Let $\lhd$ be a strong inclusion on a sub-pcd-lattice  $P$ of a locale $L$. Given a locale map $f:L\to L'$  we say that $\lhd$ is
\emph{finer than $f^-\times f^-[\prec]$} if $y\prec x$ on $L'$ implies $f^-[y]\leq p\lhd p'\leq f^-[x]$, for some $p,p'\in P$.
Denote by $\cal C$  the class of mappings $f:L\to L'$ with compact and regular codomain such that $\lhd$ is finer than $f^-\times
f^-[\prec]$.

\begin{thm}\label{gammacomp}
Let $L$ be a locale, $P$ any sub-pcd-lattice of $L$, and $\lhd$ a strong inclusion on $P$.
Then, the locale mapping $\mu:L\to {\cal R}(P)$, defined, for $I\in B_d$, by $\mu^-(I)= \bigvee I$,
satisfies the following property: for every $f:L\to L'$
in $\cal C$ a unique locale mapping $g:{\cal R}(P)\to L'$ exists such that
$g\circ \mu= f$:

\[
\xymatrix{
L \ar[r]^{\mu} \ar[dr]_{ f} & {\cal R}(P)  \ar@{-->}[d]^{g} \\
 & L'
}
\]

\noindent If $\lhd$ is compatible, $\mu:L\to {\cal R}(P)$ is a compactification of $L$.
\end{thm}

\be{proof}
Standard  arguments show that $\mu^-$ satisfies the conditions on continuous maps. To prove the extension property, if $f:L\to L'$ is in $\cal C$, let $B'$ be a sub-pcd-lattice and basis of $L'$ ($B'$ always exists, cf. Section \ref{strongregularity}). For $a\in B'$, set $$g^- (a)=\{c\in P:(\exists b\in B') b\prec a, \ c\leq f^- (b)\}.$$
First we prove that $g^- (a)\in {\cal R}(P)$. By Restricted Separation, $g^- (a)$ is a set, and trivially an ideal. Let $c\in g^- (a)$. Then $c\leq f^- (b),$ for $ b\prec a$. By the interpolation property of $\prec$ on $B'$ (Lemma \ref{lemmacr}), there is $d\in B'$ such that $b\prec d\prec a$.
By hypothesis, there are $p,p'\in P$ such that $f^-(b)\leq p\lhd p'\leq f^-(d)$. Then $c\lhd p'\leq f^-(d)$; since $d\prec a$, one concludes
$p'\in g^- (a)$. This gives $g^- (a)\in {\cal R}(P)$.

Now assume $a\in B'$, $U\in \Pow(B')$ and $a\leq \bigvee U$. Then $g^-(a)\leq \bigvee_{q\in U} g^- (q)$: let $c\in g^- (a)$, i.e., assume
there is $b\in B'$ such that $c\leq f^- (b), \ b\prec a$. By Lemma \ref{lemmacr}, either $b=0$, whence $c=0$, or there are finite subsets $\{q_1,...,q_n\}$, $\{q'_1,...,q'_n\}$ and $\{q''_1,...,q''_n\}$ of $B'$,  with $q'_i\prec q''_i\prec q_i$, $q_i\in U$, and such that $b\leq \vee q_i' \prec \vee q''_i\prec
\bigvee U.$
By hypothesis, $p_i,p_i'$ exist such that $f^-(q'_i)\leq p_i\lhd p_i'\leq f^-(q_i'')$. Thus, as $q_i''\prec q_i$,  $p_i'\in g^-(q_i)$.
Since $c\leq f^- (b)\leq f^- (\vee q'_i)\leq \vee p_i'$, one has $c\in \bigvee_{q\in U} g^- (q)$. Therefore $g^-
(a)\leq \bigvee_{q\in U} g^- (q)$.

By the properties of $\prec$, and the fact that $f^-$ defines a continuous map, it immediately follows that $g^-$ also satisfies conditions 1, 2 on continuous maps (since $B'$ is a pcd-lattice, to verify condition 2 it is enough to show that $g^-(a) \wedge g^-(b) = g^-(a \wedge b)$, as $g^-(a \wedge b) = \bigvee \{g^-(c): c\in B', c\leq  a, c\leq b\}$, for all $a,b\in B'$). So $g:{\cal R}(P)\to L'$ is a continuous map.

To show that $f^- (a)=\mu^- [g^- (a)]$, for all $a\in B'$, note first that, for $I\in {\cal R}(P)$, $\mu^-(I)=\bigvee I$;
  to prove then  $\bigvee g^- (a)= f^- (a)$, let $c\in g^- (a)$, i.e., assume
there is $b\in B'$ such that $ c\leq f^- (b), \ b\prec a$. This implies $c\leq f^- (a)$, so that $\bigvee g^- (a)\leq f^- (a)$.
Conversely, since $a=\bigvee \{b\in B':b\prec a\}$, $f^- (a)=f^-(\bigvee \{b\in B':b\prec a \})=\bigvee_{b\prec a, b\in B'} f^- (b)$. Let
$b\prec a$; by Lemma \ref{lemmacr}, there is $d\in B'$ such that $b\prec d\prec a$. By hypothesis, there are $p,p'\in P$ such that
$f^-(b)\leq p\lhd p'\leq f^-(d)$. Since $d\prec a$, one gets $p'\in g^- (a)$. Thus $f^-(b)\leq \bigvee g^- (a)$ for all $b\prec a$, so
that $f^- (a)\leq \bigvee g^- (a)$. Finally, uniqueness of $g$  follows from Lemma \ref{denseoneone}, $iii$, since $\mu$ is trivially dense.

If $\lhd$ is compatible,  recalling  Lemma \ref{representation},
to prove that $\mu$
is a compactification we only have to show that
$\mu^-$ is onto: by compatibility of  $\lhd$, $L$ is set-generated by $P$, and  given $a\in P$,
the set $\{b\in P: b\lhd a\}$ is a round ideal of $P$ whose join is $a$.
\en{proof}

\medskip

Now let $\F\equiv \{f_i\}_{i\in I}$ be a (possibly empty) family of locale mappings, $f_i:L\to L_i$, with $L_i$ compact and regular, and
let $S$ be a subset of $L$. Set $S_\F= S \cup \{f_i^- (b):i\in I, b\in B_i\}$,  where $B_i$ is a set of generators and sub-pcd-lattice of
$L_i$. With $S_{\mathcal F}^*$ we shall abbreviate $(S_{\mathcal F})^*$.

\begin{lem}\label{siF}
The relation $\lhd_\F \subseteq S^*_\F \times S^*_\F$
defined  inductively as the least relation containing the set $\{(f_i^-(b), f^-_i(a)): i\in I, a,b\in
B_i,  b\prec a\}$, and closed under conditions from 1  to 5  on strong inclusions,  is a strong inclusion on $S^*_\F$.
\end{lem}

\be{proof} The set of pairs $(f_i^-(b),$ $f^-_i(a))$ is a relation satisfying the hypotheses of Lemma \ref{idsi} (as
$f_i^-$ preserves $\prec$, and since, by Lemma \ref{lemmacr}, $\prec$ interpolates on $B_i$).
\en{proof}
\medskip

\noindent Starting from $S=\emptyset$, one obtains the least strong inclusion $\lhd_\F$ containing the indicated set, on the least
(sub)structure on which it may be considered, i.e., $\emptyset^*_\F$.

\begin{prop}
\label{extension} Let $L$ be a locale, $S$ any subset of $L$, $\F$ a family of locale mappings as above. Let $\lhd$ be a strong inclusion on $S_\F^*$ containing $\lhd_\F$. Then the locale mapping $\mu:L\to {\cal R}(S_\F^*,\lhd)$, $\mu^-(I)= \bigvee I$, is such that each $f_i$ factors uniquely through $\mu$: for all $i$ a unique $g_i:{\cal R}(S^*_\F,\lhd)\to L_i$ exists such that $g_i\circ \mu=f_i$.
\end{prop}

\be{proof}  By Lemma \ref{lemmacr}, if $b\prec a$ on $L_i$, one has   $b\leq b'\prec a'\leq a$, for $a',b'\in B_i$. Then,
by construction of $S_\F^*$ and $\lhd_\F$, each $f_i$ is such that $\lhd$ is finer than $f_i^-\times f_i^-[\prec]$, and hence, by
Theorem \ref{gammacomp}, each $f_i$ factors uniquely through $\mu$.
\en{proof}

\medskip

The compactification to be described next was defined in \cite{CuSCAlexandroff} using topos-valid, but impredicative, comprehension principles. There I claimed that co-inductive definitions would probably have allowed for a predicative construction\footnote{At the time of writing \cite{CuSCAlexandroff}, \cite[Theorem 13.2.3]{AR2010}, on the existence of the largest set coinductively defined by a given set-sized inductive definition, had not been (proved or) published yet.}. The following shows that this is indeed the case.

In general, the strong inclusion $\lhd_\F$ one obtains by Lemma \ref{siF}  need not be compatible, so that
$\mu:L\to {\cal R}(B_\F^*, \lhd_\F)$ need not be a compactification. For obtaining a compatible strong inclusion, we shall have to enlarge $\lhd_\F$.

\begin{thm}
\label{gammaext} Let $L$ be a strongly regular locale,  and let $\F\equiv \{f_i\}_{i\in I}$ be a (possibly empty) family of locale
mappings $f_i:L\to L_i$, with $L_i$ compact and regular for all $i\in I$. Let $B$ be a strongly regular basis for $L$, and $B_i$  any  basis for  $L_i$.
Define  $B_\F = B \cup \bigcup_{i\in I} \{f_i^- (b_i):b_i\in B^*_i\}$. On  $B_\F^*$ let $\lhd_\prec$ be the largest interpolative
relation contained in $\prec\cap B_\F^*\times B_\F^*$. Then the following hold:

\begin{enumerate}

    \item[$i.$] $\lhd_\prec$ is a set and a strong inclusion on  $B_\F^*$  compatible with $L$;

    \item[$ii.$] $\mu:L\to \gamma L = {\cal R}(B_\F^*,\lhd_\prec)$ is a compactification of $L$ such that for all $f_i$ there is a unique $g_i:\gamma L\to L_i$ such that  $g_i\circ \mu= f_i$:

\[
\xymatrix{
L \ar[r]^{\mu} \ar[dr]_{ f_i} & \gamma L  \ar@{-->}[d]^{g_i} \\
 & L_i
}
\]

\end{enumerate}

\end{thm}

\be{proof}
That $\lhd_\prec$ is a set follows from proposition \ref{largestinterpolativesubrelation}, and that it is a strong inclusion is proved after the definition of strong inclusion.

We show that  $\mu:L\to \gamma L$ is a compactification of $L$ by proving  that $\lhd_\prec$ is compatible.
Since $B$ is a basis, for all $a\in B_\F^*$, $a=\bigvee_L \{b\in B: b \leq a\}$.
On the other hand, for all $b\in B$, $b=\bigvee_L \{c\in B: c \prec_0' b\}$, because  $B$ is a strongly regular basis of $L$.
One has $c \prec_0' b$ implies $c\lhd_\prec b$, since, by $B\sub B_\F^*$, $ \prec_0' \sub \prec \cap (B\times B)\sub \prec\cap  (B_\F^*\times B_\F^*)$, and $\prec_0'$ is interpolative by definition.  Then, $b=\bigvee_L \{c\in B: c \lhd_\prec b \}$, whence $a=\bigvee_L \{b\in B_\F^*: b \lhd_\prec a\}$. This proves that $\lhd_\prec$ is compatible, so that $\mu$ is a compactification.

Since $\lhd_\F$ is interpolative and contained in $\prec\cap (B_\F^*\times B_\F^*)$, one has that $\lhd_\prec$ contains
$\lhd_\F$. Thus, by Proposition \ref{extension}, $ii$ follows.
\en{proof}

\medskip

We showed that a compatible strong inclusion on a sub-pcd-lattice $P$ of a locale $L$ determines a compactification of $L$. The converse also holds. To prove this we shall need the following lemma.

\begin{lem}\label{si_explicitdefinition}
Let $\F=\{f\}$, $f:L\to M$ with $M$ compact and regular, $S$ any subset of $L$, and let $B_M$ be a base and sub-pcd-lattice of $M$. If $f^-[\cdot]:M\to L$ is $*-$preserving, then for all $x,y\in S^*_f$,  $x\lhd_f y \iff x\leq f^-(b), f^-(a)\leq y$ for some $a,b\in B_M$, with $b\prec a$.
\end{lem}

\be{proof} On $P=S^*_f$, $\lhd_f$ is defined inductively as in Lemma \ref{idsi}, as the least relation containing the set $R=\{(f^-(b), f^-(a)):  a,b\in
B_M,  b\prec a\}$, and closed under conditions from 1  to 5  on strong inclusions. I.e., $\lhd_f= I(\Phi_R)$, with $\Phi_R$ as in Lemma \ref{idsi}.
We prove the claim by induction. Let $T\subseteq P\times P$, $T= \{(x,y)\in P\times P: (\exists a\in B_M)(\exists b\in B_M)b\prec a \ \&\ x\leq f^-(b)\ \&\ f^-(a)\leq y\}$. We prove that $T$ is $\Phi_R-$closed, so that, since $\lhd_f$ is the least  $\Phi_R-$closed class, we may conclude $\lhd_f\subseteq T$.
Given $(f^-(b), f^-(a))$ with $a,b\in
B_M,  b\prec a$, we have to prove that $(\exists c\in B_M)(\exists d\in B_M)d\prec c \ \&\ f^-(b)\leq f^-(d)\ \&\ f^-(c)\leq f^-(a)$. One simply takes $d=b,c=a$. The remaining cases are as simple, we only prove two of them: assume $x\leq f^-(b'), f^-(a')\leq a$ with $b'\prec a'$, and $y\leq f^-(b''), f^-(a'')\leq a$ with $b''\prec a''$. Then, $x\vee y\leq f^-(b')\vee f^-(b'')=f^-(b'\vee b''), f^-(a'\vee a'')=f^-(a')\vee f^-(a'')\leq a$, with $b'\vee b''\prec a'\vee a''$, so that $(x\vee y,a)\in T$. If instead $a\leq f^-(b'), f^-(a')\leq b$, with $b'\prec a'$, then $b^*\leq f^-(a')^*=f^-(a'^*)$ and $f^-(b')^*=f^-(b'^*)\leq a^*$, with $a'^* \prec  b'^*$, since $f^-$ is $*-$preserving.
\en{proof}

\begin{thm}\label{isom}
Let $k:L\to kL$ be a compactification of a locale $L$. Then a sub-pcd-lattice $P_k$ of $L$ and a compatible strong inclusion $\lhd_k$ on $P_k$ exist such that ${\cal R}(P_k,\lhd_k)$ is isomorphic with $kL$.
\end{thm}

\be{proof}
Let $P_k=\emptyset^*_{\{k\}}$, and $\lhd_k$ be the strong inclusion given by Lemma \ref{siF}, where $\F= \{k\}$, i.e., the least relation on $\emptyset^*_{\{k\}}$ containing the set $\{(k^-(b), k^-(a)): a,b\in B_{kL},  b\prec a\}$, and closed under conditions from 1  to 5  on
strong inclusions ($B_{kL}$ sub-pcd-lattice and basis of $kL$).
Compatibility of  $\lhd_k$: for all $a\in P_k$, $a=\bigvee_{b\in U} k^-(b)$ for some $U\in \Pow(B_{kL})$, as $k^-$ is onto.
For $b\in U$, one has $b=\bigvee_{c\in B_{kL}}\{c: c\prec b\}$.
Thus, $a=\bigvee_{b\in U}\bigvee_{c\in B_{kL}}\{k^-(c):c\prec b\}$. Since $k^-(c)\lhd_k k^-(b)\leq a$, one concludes that $\lhd_k$ is compatible.

It remains  to prove that $kL$ is isomorphic with ${\cal R}(P_k,\lhd_k)$. By Theorem \ref{gammacomp}, the mapping $\mu:L\to {\cal R}(P_k)$, defined by $\mu^-(I)= \bigvee I$, is a compactification of $L$ that satisfies: for every $f:L\to L'$ in $\cal C_k$ a unique mapping $g:{\cal R}(P_k)\to L'$ exists such that
$g\circ \mu= f$, with $\cal C_k$  the class of mappings $f:L\to L'$ with compact and regular codomain such that $\lhd_k$ is finer than $f^-\times
f^-[\prec]$. Now, compactification $k$ is in $\cal C_k$. Indeed, if $y\prec x$ in $kL$, by Lemma \ref{lemmacr} there are $q,q'\in  B_{kL}$ with $y\leq q\prec q'\prec x$. Therefore, $k^-[y]\leq k^-(q)\prec k^-(q')\prec k^-[x]$, with $k^-(q), k^-(q')\in P_k$. We then have a unique $g:{\cal R}(P_k)\to kL$ such that $g\circ \mu= k$, defined by, for $a\in B_{kL}$,
$$g^- (a)=\{c\in P_k:(\exists b\in B_{kL}) b\prec a, \ c\leq k^- (b)\}.$$
We prove that $g$ provides the required isomorphism. Since $k$ is dense, $g$ is dense too (as $k=g\circ \mu$). By Lemma \ref{denseoneone}, then $g^-[\cdot]$ in one-one. To conclude thus we have to prove that $g^-$ is onto: if $I\in  {\cal R}(P_k,\lhd_k)$ there is $U\in \Pow(B_{kL})$ with $\bigvee_{b\in U} g^-(b)=I$. It is enough to show that for $a\in P_k$ there is $U\in \Pow(B_{kL})$ with $\bigvee_{b\in U} g^-(b)=\Downarrow a$. Let $U=\{a'\in B_{kL}: k^-(a')\leq a\}$.
If $d\in \Downarrow a$, then $d\lhd_k a$, so that by the previous Lemma, $d\leq k^-(b'), k^-(a')\leq a$ for some $a',b'\in B_{kL}$, $b'\prec a'$ ($k^-$ is $*-$preserving by Lemma \ref{denseoneone}, $i.$). Thus, $\Downarrow a\subseteq \bigcup_{a'\in U} g^-(a')=\bigcup_{a'\in U}\{c\in P_k:(\exists b'\in B_{kL}) b'\prec a', \ c\leq k^- (b')\}$, that gives $\Downarrow a\leq \bigvee_{a'\in U} g^-(a')$. On the other hand, for all $a'\in U$, $g^-(a')\leq \Downarrow a$ again by the previous Lemma. This proves that $g^-$ is onto, as wished.
\en{proof}

Note that the class $\cal C_k$, of continuous maps $f:L\to L'$ with compact and regular codomain such that $\lhd_k$ is finer than $f^-\times f^-[\prec]$ (recall that $\lhd_k$ is
finer than $f^-\times f^-[\prec]$ if $y\prec x$ on $L'$ implies $f^-[y]\leq p\lhd_k p'\leq f^-[x]$, for some $p,p'\in P_k$), is completely determined by the compactification $k:L\to kL$, i.e., it does not depend on the choice of the base $B_{kL}$ of $kL$ for the construction of $(P_k,\lhd_k)$. Indeed, assume $B'_{kL}$ is a different sub-pcd-lattice and basis of $kL$, and let $P_k'=\emptyset^*_{\{k\}}$ and $\lhd'_k$ be the sub-pcd-lattice of $L$ and strong inclusion constructed from $B'_{kL}$. We prove that if $\lhd'_k$ is finer than $f^-\times f^-[\prec]$, then also $\lhd_k$ is finer than $f^-\times f^-[\prec]$:
let $y\prec x$ on $L'$, so that  $f^-[y]\leq p\lhd'_k p'\leq f^-[x]$ for $p,p'\in P_k'$.
By Lemma \ref{si_explicitdefinition}, if $p\lhd'_k p'$ then $a',b'\in B'_{kL}$ exist with $b'\prec a'$, $p\leq k^-(b')$, $k^-(a')\leq p'$, as $k^-$ is $*-$preserving. Since $kL$ is compact, by Lemma \ref{lemmacr}, $a,b$ exist in $B_{kL}$ such that  $b'\prec b\prec a\prec a'$. Therefore, $p\leq k^-(b)=q$, $q'=k^-(a)\leq p'$, whence $q\lhd_k q'$, and $f^-[y]\leq q\lhd_k q'\leq f^-[x]$.
\smallskip

Banaschewski \cite{B90} proves (in our notation) that every compactification  $k:L\to kL$ is associated with a compatible strong inclusion $\lhd_k$ on $L$ determining an isomorphic compactification ${\cal R}(L,\lhd_k)$, by defining, for $x,y\in L$, $x \lhd_k y\iff \bigvee\{u\in kL:k^-(u)=x\}\prec \bigvee\{v\in kL:k^-(v)=y\}$. As already expounded, Banaschewski's  construction is not viable in  the present constructive setting. Theorem \ref{isom} provides then a constructive version of Banaschewski's result.

The construction of ${\cal R}(P_k,\lhd_k)$ in the previous theorem, aside from its constructive character, allows us to characterize a compactification in terms of its ability  to extend continuous mappings with compact and regular codomain.
Recall indeed that the compactifications of a locale $L$ are ordered as follows: if $k:L\to kL$ and $k':L\to k'L$ are compactifications of $L$, $k\leq k'$ if and only if $k=h\circ k'$, for a continuous $h:k'L\to kL$:

\vspace{1cm}\hspace{2cm}
 $k\leq k'\hspace{1cm} \iff$

\vspace{-1.5cm}

\begin{align*}\hspace{2cm} \xymatrix{
&L \ar[ld]^{k}\ar[rd]_{k'} & \\
kL && k'L \ar[ll]^h}
\end{align*}

\noindent The following is the announced characterization.

\begin{cor}\label{characterization}
Let  $k:L\to kL$ be a compactification of a strongly regular locale $L$, and $\cal C_k$ the associated class of continuous maps $f:L\to L'$ with compact and regular codomain such that $\lhd_k$ is finer than $f^-\times f^-[\prec]$. Then, for every $f\in \cal C_k$ a unique continuous $g:kL\to L'$ exists such that $g\circ k= f$, and $k:L\to kL$ is the minimal compactification with this property.
\end{cor}
\be{proof}
The first part of this corollary follows directly from Theorems \ref{gammacomp} and \ref{isom}. For the second, assume $k':L\to k'L$ has the same property with respect to the class $\cal C_k$. Then, since $k:L\to kL$ is in $\cal C_k$ (cf. proof of Theorem \ref{isom}), there is a (unique) $g:k'L\to kL$ with $g\circ k'= k$, so that $k\leq k'$.
\en{proof}

\smallskip

Observe that, given a strong inclusion $\lhd$ on a sub-pcd-lattice $P$ of a locale $L$, the composite  $\overline \lhd\equiv \leq \circ \lhd \circ
\leq$ on $L$ yields a class relation satisfying the axioms of a  strong inclusion, the least one on $L$ extending
$\lhd$. One then has that $\lhd$ on $P$ is finer than $f^-\times f^-[\prec]$ iff   $f^-\times f^-[\prec]\subseteq \overline \lhd$. So the extension property in the previous corollary may be expressed as follows: every continuous $f:L\to L'$  with compact and regular codomain such that $f^-\times f^-[\prec]\subseteq \overline \lhd_k$ has a unique continuous extension $g:kL\to L'$ to $kL$.
In \cite{CuSCAlexandroff} we showed that Alexandroff compactification $\alpha L$ of a locally compact regular locale $(L,B)$ can be defined by constructing inductively the least strong inclusion $\lhd_{<<}$ on $B^*$ that contains the restriction of the familiar way-below relation to $B^*\times B^*$. Once the concepts of local compactness and way-below relation are treated as indicated in \cite{CuExistenceSC}, the given proof can be formulated in CST.
So by the above Corollary, Alexandroff compactification $\alpha L$ can be characterized as
the minimal compactification that allows for the extension of all those continuous mappings from $L$ to a compact and regular codomain whose inverse images `send the well-inside relation into the least strong inclusion containing the way-below relation'.
\smallskip

Finally, we observe that for compactifications $k_1:L\to k_1 L,k_2:L\to k_2 L$ one has  $$k_1\leq k_2  \iff \cal C_{k_1}\subseteq \cal C_{k_2}$$
Indeed, assume $k_1\leq k_2$, i.e., assume there is $h:k_2 L \to k_1 L$ such that $k_1=h\circ k_2$. Let $f:L\to L'$ in $\cal C_{k_1}$, so that, if $y\prec x$ in $L'$, there are $p,p'\in P_{k_1}$ with $f^-[y]\leq p\lhd_1 p'\leq f^-[x]$. To conclude that $f\in \cal C_{k_2}$, we want $q,q'\in P_{k_2}$ such that $f^-[y]\leq q\lhd_2 q'\leq f^-[x]$. By $p\lhd_1 p'$, recalling Lemma \ref{si_explicitdefinition}, we have $a,b\in B_{k_1L}, b\prec a$ with $p\leq k_1^-(b)=k_2^-[h^-(b)],\  k_2^-[h^-(a)]=k_1^-(a)\leq p'$ ($B_{k_iL}$ sub-pcd-lattices and bases for $k_iL$). Since $k_2L$ is compact and regular, and since $h^-(b)\prec h^-(a)$, by Lemma \ref{lemmacr} we have $a',b'\in B_{k_2L}$ with $h^-(b)\prec b'\prec a'\prec h^-(a)$, so that $f^-[y]\leq p\leq k_1^-(b)=k_2^-[h^-(b)]\leq k_2^-(b')\lhd_2 k_2^-(a')\leq k_1^-(a)\leq p'\leq f^-[x]$. Thus we conclude with $q=k_2^-(b'), q'= k_2^-(a')\in P_{k_2}$. Conversely, if $\cal C_{k_1}\subseteq \cal C_{k_2}$, we conclude by Corollary \ref{characterization} that $k_1\leq k_2$ observing that $k_1\in C_{k_1}$.

Therefore, if the sublattice $P_k$ of $L$ and strong inclusion $\lhd_k$ are `sufficiently large' to be such that$\lhd_k$ is finer than every continuous function with compact regular co-domain, so that
$\cal C_k$  coincides with the class of all such functions,
$k:L\to kL$ is the largest compactification of $L$, i.e., the compact regular reflection of $L$ (`weak' Stone-\v{C}ech compactification). In this sense, the ability of a compactification $k$ to extend the continuous functions in $\cal C_k$ can be regarded as an approximation of the characterizing property of Stone-\v{C}ech compactification (when it exists, cf. next section).

\section{An interesting open question: the constructive compact regular reflection of a locale}
\label{openproblem}

One might be tempted to think that, choosing in Theorem \ref{gammaext} as $\F$ the family of all mappings $f : L \to L'$, with $L'$ compact regular, one would get the `weak' Stone-\v{C}ech compactification of $L$ (i.e., its compact regular reflection, `weak' as opposed to the compact completely regular reflection). Unfortunately, this is not the case, since those mappings are too many to be set-indexed (even in classical set theory).

In \cite{BP}, Banaschewski and Pultr define the compact regular reflection of a strongly regular locale $L$ considering the strongly regular ideals of $L$ with respect to the largest strong inclusion contained in the well-inside relation on $L$. Previously, Banaschewski and Mulvey \cite{BM80} used the completely regular ideals on $L$ (i.e., the strongly regular ideals with respect to the really inside relation as strong inclusion) to define the compact completely regular reflection of $L$ (`strong' Stone-\v{C}ech compactification).

None of these reflections can be defined constructively in the present sense in full generality. This is not just because the specific constructions cannot be carried out in CST (a non-trivial frame is a proper class \cite{CuExistenceSC}, and the ideals of this class do not yield a frame).
It is proved in \cite{CuExistenceSC} that the compact \emph{completely} regular reflection is independent from  CST  also when extended with (various principles including) the principle of dependent choice.
With DC the compact regular and compact completely regular reflection coincide, so that if we could define the compact regular reflection of every locale $L$ in CST, we could also define the compact completely regular reflection of every locale $L$ in CST+DC (see the Appendix for the formal version of these statements).

In \cite{CuSC} it is however proved that one may define constructively the compact completely regular reflection $\beta L$ of $L$, for every locale  $L$ such that the class of continuous mappings $Hom_{Loc}(L,[0,1])$ from $L$ to the localic real unit interval is a set (while this is always the case in classical set theory or in a topos, proving that $Hom_{Loc}(L,[0,1])$ is a set is in general non-trivial, and not the case for every locale in CST). This condition in fact characterizes the locales $L$ of which the compact completely regular reflection can be constructed.

The construction of $\beta L$ in \cite{CuSC} can be summarized as follows: if $Hom_{Loc}(L,[0,1])$ is a set, one can expand a base $B$ of $L$ to a base $B_\F$ containing counterimages of basic elements of the localic real unit interval. Since $Hom_{Loc}(L,[0,1])$ is a set, $B_\F$ is a set, too. Then, on the pcd-lattice  $B_\F^*$ one takes the completely regular ideals. The resulting frame ${\cal R}(B_\F^*)$ is compact and completely regular, and allows for the extension of each mapping in $Hom_{Loc}(L,[0,1])$ to ${\cal R}(B_\F^*)$.

One then proves that (\v{C}ech's theorem implies Stone's theorem, i.e., that) a compactification of $L$  that allows for the extension of each mapping in $Hom_{Loc}(L,[0,1])$ in fact allows for the extension of each continuous mapping in $Hom_{Loc}(L,L')$, for every $L'$ compact and completely regular. This latter result in turn is proved exploiting Tychonoff embedding theorem for exhibiting the completely regular locale $L'$ as a sublocale of a set-indexed product of the localic real unit interval.

A similar procedure does not seem to work for the compact regular reflection of a strongly regular locale. Indeed, when $Hom_{Loc}(L,[0,1])$ is a set for $L$ a strongly regular locale, by Theorem \ref{gammaext} we are able to construct a compact regular compactification that allows for the extension of all mappings in $\F=Hom_{Loc}(L,[0,1])$, but we cannot use Tychonoff embedding theorem for embedding a compact regular (not necessarily completely regular without DC) $L$ in a set-indexed product of $[0,1]$, for showing that the obtained compactification is in fact the compact regular reflection of $L$.

In \cite{BM80} Banaschewski and Mulvey also provided a definition of the compact regular reflection of a locale, but remarked that a concrete description of the compact regular reflection of a strongly regular locale using strongly regular ideals would have been desirable. As said such a description, when not appealing to strongly non constructive set existence principles, cannot be derived in full generality, but the question whether a restricted version of it similar to that obtained for the compact completely regular reflection in \cite{CuSC} can be constructed is open.
\medskip

\noindent \textbf{Open Problem:} Can the  compact regular reflection of a locale $L$  be defined in CST in non-trivial cases? Can the locales of which the compact regular reflection exists in CST be characterized?

\medskip

A related open problem is the following. In \cite{CuSCAlexandroff} we proved that, when $Hom_{Loc}(L,$ $[0,1])$ is a set, the sub-pcd-lattice $B_\F^*$ of $L$ used for constructing the compact completely regular reflection  is sufficiently large for fully replacing $L$ in the construction of every compactification, in the sense that every
compactification of a frame $L$ can be obtained as the frame of round ideals over $B_\F^*$, for a certain inductively defined strong inclusion on $B_\F^*$. A similar result for compact regular compactifications would be desirable.

\section {Appendix: Constructive Set Theory. Inductive and co-inductive definitions}\label{CST}

As pointed out before, we used CST as a collective name for Aczel-Myhill formal systems for constructive set theory.
In this appendix the specific formal system in which we have been working in is described to make the article self-contained. The reader may consult \cite{A86,AR,AR2010} for a thorough introduction to the subject.
A core formal system for CST is the choice-free  \emph{Constructive Zermelo-Fraenkel Set Theory} (CZF).
This system is often extended by principles (described below) ensuring that certain inductively and co-inductively defined classes are sets.
Note that CZF, extended or not by those principles,  is a subtheory of classical set theory. As already pointed out, in contrast to ZF, CZF does not have the impredicative unrestricted Separation scheme and the Powerset axiom.

The language of CZF is the same as that of Zermelo-Fraenkel Set Theory, ZF, with $\in$ as the only non-logical symbol. Beside the rules and axioms of a standard calculus for intuitionistic predicate logic with equality (e.g., \cite{TvD}), CZF has the following axioms and axiom schemes:

\begin{enumerate}
  \item Extensionality: $\forall a\forall b (\forall y (y \in a \leftrightarrow y \in b) \rightarrow a = b)$.
  \item Pair: $\forall a\forall b\exists x\forall y (y \in x\leftrightarrow y = a \vee y = b).$
  \item Union: $\forall a\exists  x\forall y (y \in x \leftrightarrow  (\exists z \in a)(y \in z))$.
  \item Restricted Separation scheme: \smallskip

\qquad \qquad $\forall a\exists  x\forall y (y \in x \leftrightarrow  y\in a\& \phi(y)),$
 \smallskip

 \noindent for $\phi$ a restricted formula. A formula $\phi$ is \emph{restricted} if the quantifiers that occur in it are of the form $\forall x\in b$,  $\exists x\in c$.

  \item Subset Collection scheme:
  \begin{eqnarray*}&&\forall a\forall b\exists c\forall u ((\forall x \in a)(\exists y \in b)\phi(x, y, u)\; \rightarrow\\
  &&  (\exists d \in c)((\forall x \in a)(\exists y \in d)\phi (x, y, u)
\& (\forall y \in d)(\exists x \in a)\phi (x, y, u))).\end{eqnarray*}

  \item Strong Collection scheme:
  \begin{eqnarray*}&&\forall a ((\forall x \in a)\exists y \phi(x, y)\;\rightarrow\\
      && \exists b ((\forall x \in  a)(\exists y \in b)\phi (x, y) \& (\forall y\in b)(\exists x \in a)\phi(x, y))).
   \end{eqnarray*}

  \item Infinity: $\exists a (\exists x\in a  \&  (\forall x\in a)(\exists y\in a)  x\in y)$.

  \item Set Induction scheme: $\forall a ((\forall x \in a)\phi(x)\rightarrow \phi(a))\rightarrow \forall a\phi(a)$.
\end{enumerate}

\noindent
We shall denote by $\CZF^-$ the system obtained from $\CZF$ by leaving out the Subset Collection scheme.
Subset Collection is perhaps the most unusual of the CZF axioms and schemes; for  this article it suffices to note that using it one proves that the class $b^a$ of functions  from a set $a$ to a set $b$ is a set, i.e., the  Exponentiation Axiom. We did not make use of Subset Collection or Exponentiation in this article.
Note that the theory obtained from $\CZF$ by adding the Law of Excluded Middle has the same theorems as ZF.

As in classical set theory, one takes advantage in this context of class notation and terminology \cite{AR,AR2010}. For example, given any set or class $X$, one has the class $\mathrm{V}=\{x\mid x=x\}$ of all sets.

As shown in this article, a major role in constructive set theory is played by inductive definitions.  An \emph{inductive definition}  is any class  $\Phi$ of pairs.
A class $A$ is \emph{$\Phi-$closed} if:
\smallskip

\centerline{$(a,X)\in \Phi$, and  $X\subseteq A$ implies $a\in A$.}
\smallskip

\noindent
The following  theorem is called \emph{the class inductive definition theorem} \cite{AR}. It says that (any extension of) the system $\CZF^-$ has the
class induction property CIP recalled in Section \ref{preliminaries}.

\be{thm}[$\CZF^-$]\label{cidt}
Given any class  $\Phi$ of ordered pairs, there exists a least $\Phi-$closed class $I(\Phi)$, the class inductively defined by $\Phi$.
\en{thm}

Even when $\Phi$ is a set, $I(\Phi)$ need not be a set in CZF. For this reason, CZF is often extended with the  \emph{Regular Extension Axiom}, REA.
\smallskip

\qquad REA: every set is a subset of a regular set.
\smallskip

\noindent A set $c$ is \emph{regular} if it is transitive, inhabited, and for any $u\in c$ and any set $R\subseteq u \times c$,  if $(\forall x \in u)(\exists y) \langle x, y\rangle \in R$, then there is a set $v\in  c$ such that
\begin{eqnarray}\label{rea}(\forall x \in u)(\exists y \in v)((x, y) \in R) &\& & (\forall y \in v)(\exists x\in u)((x, y)  \in R).\end{eqnarray}
$c$ is said to be {\em weakly regular} if in the above definition of regularity the second conjunct in (\ref{rea}) is omitted.
The weak regular extension axiom, wREA, is the statement that every set is the subset of a weakly regular set.

A class $K$ is  a \emph{bound} for $\Phi$ if, for every $(x,X)\in \Phi$, there is a set $k\in K$ and an onto mapping $f:k\to X$. The inductive definition $\Phi$ is defined to be \emph{bounded} if:

\noindent
1. $\{x\mid (x,X)\in \Phi\}$ is a set for every set $X$;

\noindent 2. $\Phi$ is bounded by a set.

The following theorem states that in the system $\mathbfr{CZF + wREA}$  the Bounded Induction Scheme (BIS) is derivable.

\be{thm}[$\mathbfr{CZF + wREA}$]\label{sidt}
If  $\Phi$ is bounded, in particular if $\Phi$ is a set, then $I(\Phi)$ is a set.
\en{thm}

Given an inductive definition $\Phi$ on a class $S$, $\Phi\subseteq S\times \Pow(S)$, a class $C\subseteq S$ is said to be $\Phi-$inclusive if $C\subseteq \Gamma_\Phi (C)$, with $\Gamma_\Phi$ the operator on subclasses associated with $\Phi$: $\Gamma_\Phi(C)=\{x \ \mid \ (\exists X)\ (x,X)\in \Phi \ \& \ X\subseteq C\}$. Aczel  showed that the class $J=\bigcup \{Y\in \Pow(S) \ \mid \ Y\subseteq \Gamma_\Phi (Y)\}$ is the largest $\Phi-$inclusive subclass of $S$, the class co-inductively defined by $\Phi$, denoted by $C(\Phi)$. This result is proven in the system CZF$^-$ + RRS, where the Relation Reflection Scheme RRS is the following axiom scheme.
\smallskip

\noindent Relation Reflection Scheme, RRS:

\begin{itemize}
\item[]
  For classes $S,R$ with
$R \subseteq S \times S$, if $a \in S$ and $\forall x \in S\ \exists y \in S\ R(x,y)$ then there is a set $S_0 \subseteq  S$ such that
      $a \in S_0$ and $\forall x \in S_0\ \exists y \in S_0\ R(x,y)$.
\end{itemize}

\noindent This scheme can be regarded as a weakening of  the Relativized Dependent Choices Axiom, RDC. By contrast with RDC, RRS is valid in all topological models (all cHa-valued models). Note also that RRS is a theorem of ZF (see \cite{A2008} for a proof of these facts).
\smallskip

\noindent The following strengthening of RRS and REA is used to show that $J$ is a set  when $S,\Phi$ are sets. A regular set $A$ is \emph{strongly regular} if it is closed under the union
operation,  i.e., if $\forall x \in A\ \cup x \in A$. Let $A$ be a strongly regular set. $A$ is  defined  to be \emph{RRS-strongly regular} if also, for all sets $A' \subseteq A$ and $R\subseteq  A'\times A'$, if $a_0 \in A'$ and $\forall x \in A'\ \exists y \in A'\ xRy$ then there is $A_0 \in A$ such that $a_0 \in A_0 \subseteq A'$ and $\forall x \in  A_0\ \exists y \in A_0\ xRy$.
\medskip

\noindent
RRS-$\bigcup$REA: Every set is the subset of a RRS-strongly regular set.

\medskip

\noindent The following theorem is proved in \cite{AR2010}.

\be{thm}[CZF+ RRS-$\bigcup$REA]\label{scdt}
If  $S$ and $\Phi\subseteq S\times \Pow(S)$ are sets,   the largest $\Phi-$inclusive subclass of $S$, i.e., the class $C(\Phi)\equiv J=\bigcup \{Y\in \Pow(S) \ \mid \ Y\subseteq \Gamma_\Phi (Y)\}$ co-inductively defined by $\Phi$, is a set.
\en{thm}

\noindent This is the Set Coinduction Scheme recalled in Section \ref{preliminaries}. \noindent CZF+ RRS-$\bigcup$REA is then the formal system for CST adopted in this article.

We conclude this Appendix describing more formally the independence - mentioned in Section \ref{openproblem} - of the weak  Stone-\v{C}ech compactification of locales in its full form,  from the formal system we have adopted.
In \cite{CuExistenceSC}, the compact completely regular reflection of a Boolean locale is proved to be independent from the formal system CZF, and from every extension of CZF that is consistent with the Generalized Uniformity Principle (GUP), as e.g. CZF+REA+RDC.  As pointed out to me by M. Rathjen, using realizability \cite{Rathjen05}, one shows that one such system is also CZF+$\bigcup$REA+ RDC. This system  extends the choice-free  system  CZF+RRS-$\bigcup$REA, described above, and adopted in this article as formal system for CST. So, if the compact regular reflection of a Boolean locale were derivable in CZF+RRS-$\bigcup$REA, it would also be derivable in CZF+$\bigcup$REA+RDC. However, in the presence of DC (and a fortiori of RDC), a locale is  compact regular iff it is compact completely regular, so that the compact regular and compact completely regular reflection of a locale coincide. 

\bibliographystyle{plain}

\medskip

\end{document}